\documentclass[11pt, final, conference, letterpaper, oneside, onecolumn]{IEEEtran}		

\usepackage[cmex10]{amsmath}
\usepackage{amsthm}
\usepackage{array}
\usepackage{mdwmath}
\usepackage{mdwtab}
\usepackage{cite}
\usepackage{graphicx}
\usepackage{graphics}

\theoremstyle{definition}

\newtheorem{Def}{{\mbox{$\;\;\;\;\;\,$}}Definition}[section]
\newtheorem{Th}{{\mbox{$\;\;\;\;\;\,$}}Theorem}[section]

\pdfoutput=1

\begin{document}

\title{Nonautonomous Food-Limited Fishery Model With Adaptive Harvesting}

\author{
   L. V. Idels
   \thanks{Research supported by a grant from Vancouver Island University (VIU)}.
      \\
   Department of Mathematics  \\
   Vancouver Island University  \\
   900 Fifth St. Nanaimo BC, V9S5S5  \\
   email:  {\tt{lev.idels@viu.ca}} \\\\
   A. Khokhlov\\
 Department of Mathematics \\   Chelyabinsk South Ural State University\\
   80a Lenin Ave, Chelyabinsk  454080, \\Russia \\ email:
   {\tt{Artur.Khokhlov@rocketsoftware.com}}\\\\
   Corresponding Author:  Lev V. Idels {\vspace{-1mm}}}

\maketitle \vfill{}
\newpage
\begin{abstract}
We will introduce the biological motivation of the $\gamma$-
food-limited model with variable parameters. New criteria are
established for the existence and global stability of positive
periodic solutions. To prove the existence of steady-state
solutions, we used the upper-lower solution method where the
existence of at least one positive periodic solution is obtained by
constructing a pair of upper and lower solutions and application of
the Friedreichs Theorem. Numerical simulations illustrate effects of periodic variation in
the values of the basic biological and environmental parameters and how the adaptive harvesting strategies
affect fishing stocks.
\end{abstract}
{\bf Keywords}-Fishery models, Nonautonomous Differential Equations,
Harvesting Strategies, Food-Limited Model, Adaptive Harvesting,
Periodic Solutions, Stability.\\
Math Subject Classifications: 34C25, 34D23, 92B05

\section{Time-Varying Fishery Models: Biological motivation}

Most populations experience regular or recurring fluctuations in
biological and environmental factors which affect demographic
parameters \cite{Bo1}--\cite{Bo3}, \cite{Cla}, \cite{Cush1},
\cite{Cush2}, \cite{Hs},  and mathematical models cannot ignore
for example, year-to-year changes in weather, the global climate
variability, habitat destruction and exploitation, the expanding
food surplus, and other factors that affect the population growth
\cite{Bo1}--\cite{Bra2}, \cite{Cush2}, \cite{Hut},
\cite{La}--\cite{Ro} and \cite{Vlad}. In most models of population dynamics, increases in population due to birth are assumed to be time-independent,
but many species reproduce only during a single period of the year. There are several biological parameters that can vary
seasonally, including some cyclical changes in control parameters.
For example, in  temperature or polar zones growth frequently slows down, or even
ceases in winter.
Careful analysis in \cite{Me} shows that there might be a relationship between
asymptotic recruitment and bottom temperature, i.e., stocks located in warmer waters had lower asymptotic recruitment.

Consider the following autonomous model for the harvested population with size
$N(t)$ at time $t$ :

\begin{equation}\label{Eq0}
\frac{1}{N} \frac{dN}{dt} = G(N)-E
\end{equation}
where
\begin{equation}\label{pol}
G(0) = r \,, \qquad G(K) = 0 .\nonumber
\end{equation}
We assume that $G(N)$ is strictly
decreasing, and define the intrinsic growth rate $r>0$ for $N \approx 0$, the
carrying capacity $K>0$, and the effort
function $E$.\\
The linearity assumption in the logistic model is
violated for nearly all populations, e.g. for a food-limited population Smith \cite{SM} (see \cite{Tso})
reported a snag in the classical logistic model, i.e. it did not fit
experimental data, and suggested a modification of the logistic
equation
\begin{equation}\label{Sm}
\frac{1}{N}\frac{dN}{dt}=G(N|r,K,\beta)-E,
\end{equation}
where $G(N|r,K,\beta)=r(1-N/K)(1+\beta N/K)^{-1}$, and Smith
\cite{SM} called the coefficient $(1+\beta N/K)^{-1}$ the delaying
factor.\\
Let $x(t)=\frac{N(t)}{K}$, then equation \eqref{Sm} has the form
\begin{equation}\label{Sm0}
\frac{dx}{dt}=rx\frac{1-x}{1+\beta x}-Ex \nonumber
\end{equation}
To take into consideration a crowding factor , we
introduce a new function
\begin{equation}\label{Eq1}
G(x|r,\beta,\gamma) = r(1-x^{\gamma})(1+\beta x^{\gamma})^{-1},
\end{equation}
where $\gamma > 0$.
Then equation \eqref{Eq0} takes a form
\begin{equation}\label{Smg}
\frac{dx}{dt}=rx\frac{1-x^\gamma}{1+\beta x^\gamma}-Ex.\nonumber
\end{equation}
It is clear, if $\beta=0$, the function $G(x)$ in \eqref{Eq1} is
concave and $G(x) \rightarrow -\infty$ as $x \rightarrow \infty$. If
$ \beta>0$, then $G(x)$ is sigmoidal with $- \frac{r}{\beta} < G(x)
\leq r$.\\

{\bf {Remark.}} Let $G_{1}(x)=x(1-x^{\gamma})(1+\beta
x^{\gamma})^{-1}$. Clearly (see for example Fig 1 below) function $G_{1}(x)$ has a unique maximum
on the interval $x \in (0,1)$,
\begin{figure}[ht]
    \begin{center}
    \includegraphics[width=2.5in, viewport=83 425 533 731,clip=true]{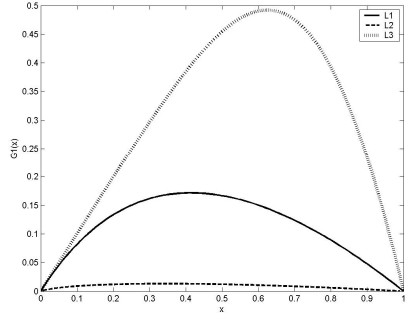}
    \label{fig1}
    \caption{Different forms of  function $G_1(x)$  (L1: $\beta=0, \gamma=1$;
    L2: $\beta=0.2, \gamma=5$; L3: $\beta=4, \gamma=0.5$) }
    \end{center}
  \end{figure}
and the maximum value and the
position of a critical point depend on some combinations of the
parameters $\gamma$ and $\beta$. If $\beta=0$ and $\gamma=1$ then
$G_{1}(x)$ is a classical symmetrical logistic function.\\
We construct a time-varying (nonautonomous) model based on equation \eqref{Sm} by allowing $r(t)$ and $K(t)$ to vary, while maintaining
constant parameters $ (\beta,\gamma)$. Similar model for unharvested population was studied in
\cite{Cush2} under the assumption that $K$ and $r$ oscillate with
small amplitudes. For an exploited marine population we introduce a varying effort function
$E=E(t)$.
\begin{equation}\label{vary}
\frac{1}{N}\frac{dN}{dt}=r(t)\left[ 1-\left( \frac{N}{K(t)}\right)
^{\gamma }\right] \left[ 1+\beta \left( \frac{N}{K(t)}\right)
^{\gamma }\right] ^{-1}-E(t),
\end{equation}
where $N|_{t=0} = N_0 >0 .$ \\In population dynamics the term
$\left(\frac{N}{K}\right)^{\gamma}$ is refereed to as the
Richards' nonlinearity \cite{Tso}. \\
 Denote
$$\overline{N}=\int_{0}^{T}N(s)ds$$
and
$$\overline{K}=\int_{0}^{T}K(s)ds ,$$
where $T > 0$ is the period of the system. Illustrative periodic functions $K(t)$  and $r(t)$ can be of
the forms:

\begin{equation}\label{rok}
r(t) = r_{0} \left[ 1 + \alpha_r \sin \frac{2\pi (t-\phi_r)}{T_r}
\right]
\end{equation}

\begin{equation}\label{K}
K(t) = K_{0} \left[ 1 + \alpha_K \sin \frac{2\pi (t-\phi_K)}{T_K}
\right]
\end{equation}
Parameters $r_{0}$, $K_{0}$, $\alpha $, and $T $  can be different for $r$ and $K$. We assume, however, that periods have a lowest common multiple,
which defines the period for the system. \\
For a canonical logistic
equation it was proven that $\overline N$ decreased as the magnitude
of variation in $K$ increases, and that $\overline N \leq \overline
K$ irrespective to $r$. However, it was shown in \cite{Cush2} that
the effect of the environmental cycles is "a very
model-dependent phenomenon".

Harvesting population models with a periodic function $E(t)$ have been studied extensively in recent
years \cite{Bra2}--\cite{Co}, \cite{Ha},
\cite{Id}-- \cite{La}. However, function $E(t)$ is not a periodic function, in fact, it is a function of many (continuous) variables which fishery managers can manipulate and this function can be defined by the fishery strategies. For example, adaptive harvesting strategies in fisheries \cite{F} rests on a combination of three elements:
 (a) deduction from prior knowledge of the ecosystem's components, (b) experience with similar ecosystems elsewhere, coupled with (c) mathematical modeling. Adaptive fishery strategies fit better with new multi-species management, with more emphasis on an "ecosystem approach" to sustainable fishery management.
The study of the dynamical behaviors of the Fox harvesting
models
\begin{equation}\label{Fox}
\frac{1}{N} \frac{dN}{dt}=r(t)\ln^{\gamma}\frac{K(t)}{N}-E(t).
\end{equation}
 in periodical environments were introduced in \cite{Bra2} and \cite{Id}.\\

The paper is organized as follows. In the next section we study qualitative
behavior of the solutions of a
harvesting model in constant environments and obtain the explicit conditions for the
existence of a unique positive solution of equation. In Section 3
for equation \eqref{vary} we will prove that it possess positive,
bounded, asymptotically stable periodic solutions. In Section 4 we
will investigate numerically the effects of periodic variation in
the values of the basic population parameters $r(t)$
and $K(t)$ and discuss adaptive  vs static harvesting strategies. \\
\section{Autonomous Model with Proportional Harvesting}
Consider equation \eqref{Smg} with proportional harvesting

\begin{equation}\label{base}
\frac{dx}{dt}=rx[ 1-x ^{\gamma }] [ 1+\beta  x ^{\gamma }] ^{-1}-Ex.
\end{equation}
If $\beta =0$ then \eqref{base} is an alternative to the logistic
fishing model with the Richards' nonlinearity
\begin{equation}\label{y}
\frac{dx}{dt}=rx( 1-x ^{\gamma }) -Ex.
\end{equation}
Equation \eqref{base} has a nontrivial  equilibrium point
$$x_{Ge}=\sqrt[\gamma]{\frac{1-E^{*}}{1+\beta E^{*}}},$$  and
equation \eqref{y} has a nonzero equilibrium point

$$x _{Le}=\sqrt[\gamma]{1-E^{*}} ,$$ where $E ^{*}=E/r $ and $0<E<r$.

For the corresponding annual equilibrium harvests
$Y_{Ge}(E)=Ex_{Ge}$ and $ Y_{Le}(E)=Ex_{Le}$, \\
 $x_{Ge}< x_{Le}$ and $Y_{Ge}\ <Y_{Le} $. Note that the maximum sustainable yield
($\max_{E}Y(E)$) exists for $E \in (0,r/q)$. \\

Let $u=\ln x$ then equations \eqref{base} has the following form $$
\frac{du}{dt}=\Phi (u)-E ,$$ where
$$
\Phi (u)=r \frac{1-e^{u \gamma}}{1+\beta e^{u \gamma}}.$$ Note
$$
\frac{d\Phi }{du}=\Phi ^{\prime }(u)=-\gamma re^{\gamma u}\frac{
1+\beta }{\left( 1+\beta e^{\gamma u}\right) ^{2}}<0
$$
and positive equilibrium
$$x^*=\sqrt[\gamma]{\frac{1-E^{*}}{1+\beta
E^{*}}}$$ of \eqref{base} is locally asymptotically stable.

\begin{Th} If in equation \eqref{base} we assume that $r,K,\beta, E $ and $ \gamma $
are all positive constants, then

$a1)$ If  $r>E$ then for every solution $x(t)$
$$
\lim_{t\rightarrow \infty }x(t)= x^{*}$$ exists.

$a2)$ If $E\equiv 0$ then $\lim_{t\rightarrow \infty }x(t)=1.$
\end{Th}
{\bf {Proof.}} The solution to equation \eqref{base} has an implicit
form
$$
E^{*}-1+(\beta E^{*} +1)x^{\gamma}=C x^{\alpha \gamma }\exp (-r\alpha
\gamma (1-E^{*})t)\
$$
 with $$\alpha =\frac{
1+E^{*}\beta  }{1+\beta },$$ $C$ is an arbitrary constant. As $
t\rightarrow \infty $
$$
x(t)\rightarrow \sqrt[\gamma ]{\frac{1-E^{*}}{1+E^{*}\beta
}}=x^{*}\text{ }.
$$

To prove the second part of the theorem, we note that if $E\equiv 0$
then the solution of \eqref{base} takes a form

$$
-1+x^{\gamma} =Cx^{\alpha\gamma}\exp (-\gamma \alpha rt)\
$$
thus $$\lim_{t\rightarrow \infty}x(t)=1 .$$

\section{Nonautonomous Model with Seasonal Harvesting}
Let us assume that all parameters in \eqref{vary} are continuous
functions and for all $t\geq 0$ \begin{equation} \label{con}
\frac{1}{N}\frac{dN}{dt}=r(t) \frac{K^{\gamma}(t)-N^{\gamma
}}{K^{\gamma}(t)+\beta N^{\gamma }}-E(t) ,
\end{equation} where $\gamma
>0$ and $\beta
>0.$

\begin{Def}
We say that a positive solution $N^{\ast}(t)$ of equation
(\ref{con}) is a global attractor or globally asymptotically stable
(GAS) if for any positive solution $N(t)$
$$
\lim_{t\rightarrow\infty}|N(t)-N^{\ast}(t)|=0.
$$
\end{Def}
Usually $N^{\ast}(t)$ is a positive equilibrium  or a positive
periodic solution of  equation (\ref{con}) if it exists. In general,
we will use the following definition.
\begin{Def}
We say that equation  (\ref{con}) is GAS, if for every two positive
solutions $N_1(t)$ and $N_2(t)$ of equation (\ref{con}) we have
$$
\lim_{t\rightarrow\infty}|N_1(t)-N_2(t)|=0.
$$
\end{Def}

Note that canonical logistic equation with variable parameters has
been well-studied, and the questions of the existence and stability
in this case are easily handled since the equation is solvable as a
Riccati equation in a closed form \cite{Bo1}--\cite{Bo3},
\cite{Cush1}, \cite{La} and \cite{Nis}.

\begin{Th}\label{teo2}
Assume that
\begin{equation} \label {eigen}
E_{1}(t)=E(t)+\frac{1}{K}\frac{dK}{dt}>0
\end{equation} and
$r(t),K(t)$ and $E(t)$ are positive functions. If $0<N(0)<K(0)$,
then every solution of equation \eqref{con} satisfies the following
inequality $0<N(t)<K(t)$ for all $t\geq 0.$
\end{Th}
{\bf {Proof.}} Let $v(t)=\ln\frac{N(t)}{K(t)},$ then equation (11)
has a form
\begin{equation} \label{mi}
\frac{dv}{dt}=r(t)\frac{1-e^{v\gamma }}{1+\beta e^{v\gamma
}}-E_{1}(t),
\end{equation}
where $v(0)<0.$ Let us prove that $v(t)<0$ for all $t\geq 0.$
Suppose there exists $t_{1}>0$ such that $v(t_{1})=0.$ Then
$$
\frac{dv}{dt} \mid_{t_{1}}=-E_{1}(t_{1})<0.
$$
Therefore in some interval ($t_{1}-a, t_{1})$ the function $v(t)$ is
decreasing and $v(t_{1})=0.$ But this is impossible, therefore
$v(t)<0$ follows by $N(t)<K(t).$ \\
To prove our next theorems for the periodic models we will use
\cite{Hart} (see also \cite{Bra2}).

\begin{Th}(Friedrichs Theorem)\label{Fr}. Suppose that $G(t,N)$ is a smooth
function with period $T$ in $t$ for every $N$. Suppose also that
there exist constants $a,b$ with $a<b$ such that $G(t,b)<0<G(t,a)$
for every $t$. Then there is a periodic solution $N_{0}(t)$ of the
differential equation $dN/dt=G(t,N)$ with period $T$ and $N(0)=c$
for some $c\in(a,b).$
\end{Th}

\begin{Th}\label{teo}
Consider equation \eqref{vary}. For all $\gamma >0,\beta \geq 0$ and
$t\geq 0$ we assume that $E(t),r(t),K(t)$ are all positive
$T-$periodic functions, and $r(t)-E_{1}(t)>0$
for all $t\geq 0.$ Then there exists a positive nonconstant periodic
solution $N_{0}(t)$ such that
$$
K(0)e^{b_{0}}<N_{0}(0)<K(0),
$$
where
$$
b_{0}<\frac{1}{\gamma }\min\lim_{t\in \lbrack 0,T]}\ln
\frac{r(t)-E_{1}(t) }{r(t)+\beta E_{1}(t)}<0.
$$\end{Th}

{\bf {Proof.}} Let $v=\ln\frac{N(t)}{K(t)},$ then equation
\eqref{vary} has a form \begin{equation}\label{op}
\frac{dv}{dt}=r(t)\frac{1-e^{v\gamma }}{1+\beta e^{v\gamma
}}-E_{1}(t)=\Phi (t,v) .\end{equation}
 Clearly, $\Phi
(t,0)=-E_{1}(t)$ $<0$ and

$$
\Phi (t,b_{0})=r(t)\frac{1-e^{b_{0}\gamma }}{1+\beta e^{b_{0}\gamma
}} -E_{1}(t)\ >0. $$ The last inequality is equivalent to the
inequality

$$
\text{\ }\frac{1}{\gamma }\ln \frac{r(t)-E_{1}(t)}{r(t)+\beta
E_{1}(t)}
>b_{0}.
$$
Therefore based on Theorem \eqref{Fr}, there exists periodic
solution $v_{0}(t)$ of equation \eqref{op}, such that
$b_{0}<v_{0}(t)<0.$ That yields the existence of the periodic
solution $N_{0}(t)$ of equation \eqref{vary} such that
$$
K(0)e^{b_{0}}<N_{0}(0)<K(0).\text{\ }
$$
\begin{Th}\label{Base} Suppose all conditions of Theorem \eqref{teo} hold.
Then there exists a unique periodic solution $N_{0}(t)>0$ of
equation \eqref{vary} such that $ 0<N_{0}(0)<K(0)$ \ and
$$
\lim_{t\rightarrow \infty }[N(t)-N_{0}(t)]=0 .$$
\end{Th}

{\bf {Proof.}} According to Theorem \eqref{teo2} the
solution  $v(t)$ of equation \eqref{mi} is negative for all $ t\geq 0.$
Firstly, let us prove that

\begin{equation}\label{ika}
\lim_{t\rightarrow \infty }\inf v(t)>-\infty .\end{equation} Assume
that

$$
\lim_{t\rightarrow \infty }v(t)=-\infty.
$$
Then
$$\lim_{t\rightarrow \infty }\inf
\frac{dv}{dt}=\lim_{t\rightarrow \infty }\inf [r(t)-E_{1}(t)]>0, $$

and $$\lim_{t\rightarrow \infty }v(t)=\lim_{t\rightarrow \infty
}\left[v(0)+\int_{0}^{t}\frac{dv}{ds}\right]=\infty .$$ We have a
contradiction. If $\lim_{t\rightarrow \infty }v(t)$ does not
exist, then there exists a sequence $\{t_{n}\}$ such that
$v(t_{n})\rightarrow -\infty $ and
\begin{equation}\label{mom}\frac{dv}{ dt}|_{t_{n}}=0.
\end{equation}
 Then equality \eqref{mom} yields
$$r(t_{n})-E_{1}(t_{n})\rightarrow 0.$$ That contradiction proves
statement \eqref{ika}. Suppose $N_{1}(t)$ and $N_{2}(t)$ are two
positive solutions of equation \eqref{mi} with $0<N_{i}(t)<K(0)$
$i=1,2$. Then $$v_{1}(t)=\ln \frac{N_{1}(t)}{K(t)}$$ and
$$v_{2}(t)=\ln \frac{N_{2}(t)}{K(t)}$$ are two solutions of equation
\eqref{mi}. Assume that
$v_{1}(t)<v_{2}(t).$ Let
$$v(t)=v_{2}(t)-v_{1}(t).$$ Then equation \eqref{mi} takes the form

\begin{equation} \label{oka}
\frac{dv}{dt}=r(t)\left[\frac{1-e^{v_{2}\gamma }}{1+\beta
e^{v_{2}\gamma }}-\frac{ 1-e^{v_{1}\gamma }}{1+\beta e^{v_{1}\gamma
}}\right].\end{equation} Application of the Mean-Value Theorem
transforms equation \eqref{oka} to \begin{equation}\label{mom}
\frac{dv}{dt}=-a(t)v(t)\end{equation} with
$$
a(t)=r(t)\left\{\frac{(1+\beta )\gamma e^{c(t)\gamma }}{[1+\beta
e^{c(t)\gamma }]^{2}}\right \}, $$ where
$$\int_{0}^{\infty}a(s)ds=\infty$$ and $$v_{1}(t)<c(t)<v_{2}(t).$$ Based on the last inequality
$$\liminf_{t\rightarrow \infty }v(t)>-\infty , $$
 $$\liminf_{t\rightarrow \infty }c(t)>-\infty .$$
Hence
$\liminf_{t\rightarrow \infty }a(t)>0,$ therefore for every
solution of \eqref{mom} we have $\lim_{t\rightarrow \infty
}v(t)=0,$ or
$ \lim_{t\rightarrow \infty }[v_{1}(t)-v_{2}(t)]=0.$
Similarly, if $v_{1}(t)> v_{2}(t)$, then $ \lim_{t\rightarrow \infty }[v_{2}(t)-v_{1}(t)]=0.$
Summing up, we conclude
$\lim_{t\rightarrow \infty }[N_{1}(t)-N_{2}(t)]=0. $

{\bf {Remark.}} If $N_{1}(t)$ and $N_{2}(t)$ are any two solutions
of equation \eqref{vary}, then statement of Theorem 3.4 is true without
the assumption that all functions $r(t)$, $K(t)$ and $E(t)$ are
$T$-periodic.\\
\newpage
\section{Numerical Experiments}
Fig 2 illustrates dynamics of the population  for a different set of parameters $\alpha$ and $\beta$.
\begin{figure}[ht]
    \begin{center}
    \includegraphics[width=2.5in, viewport=67 309 560 717,clip=true]{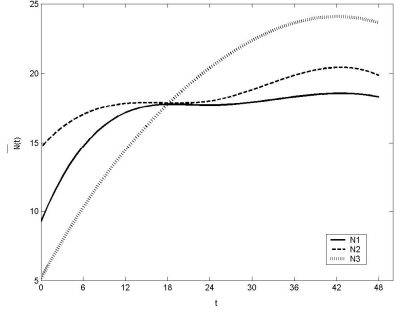}
    \label{fig2}
    \caption{Size of the population for different models: N1: $\beta=0, \gamma=1$ logistic;
    N2: $\beta=0.2, \gamma=5$; N3: $\beta=4, \gamma=0.5$ }
    \end{center}
  \end{figure}

 Fig 3  supports a well-known feature of the population models:
increase in the amplitude of the carrying capacity $K(t)$, defined by equation \eqref{K}, yields
decrease in the average size of the population, whereas,  \begin{figure}
    \begin{center}
    \includegraphics[width=2.5in, viewport=57 368 544 700,clip=true]{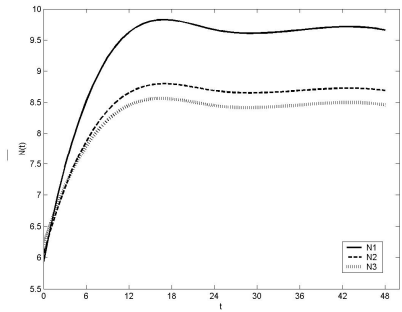}
    \label{fig3a}
  \caption{Change of environment. Population size for different values of $\alpha_{K}$
   (N1: $\alpha_{K}=0.1$, N2: $\alpha_{K}=0.5$ and N3: $\alpha_{K}=0.7$).}
    \end{center}
  \end{figure}

  Fig 4 proves that the qualitative behavior of the system is unchanged  \cite{Nis} by oscillations
in $r(t)$, defined by equation \eqref{rok}, alone.\\
{\bf {Remark.}} In all numerical experiments below we used parameters $E(t)$ and $K(t)$ that satisfy the condition \eqref{eigen}.
  \newpage
  The qualitative behavior of the system depends critically on the phase difference of the oscillations in $r(t)$ and $K(t)$.
  For example, if functions $r(t)$ and $K(t)$ have a $180^{\circ}$ phase shift, then
  the relationship between the average population size and environmental and demographic variations (see Fig 5) is
  very different from the corresponding Fig 3.
  \begin{figure}
    \begin{center}
    \includegraphics[width=2.5in, viewport=69 397 524 725,clip=true]{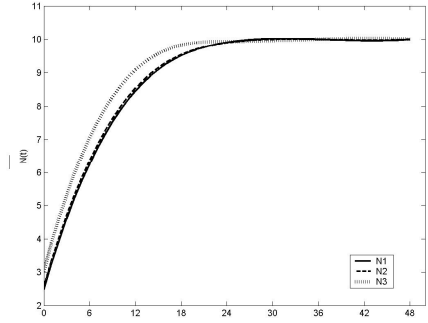}
    \label{fig4a}
  \caption{Change of demography. Population size for different values of $\alpha_{r}$
   (N1: $\alpha_{r}=0.1$, N2: $\alpha_{r}=0.5$ and N3: $\alpha_{r}=0.9$).}
    \end{center}
  \end{figure}
  \begin{figure}
    \begin{center}
    \includegraphics[width=2.5in, viewport=57 368 544 700,clip=true]{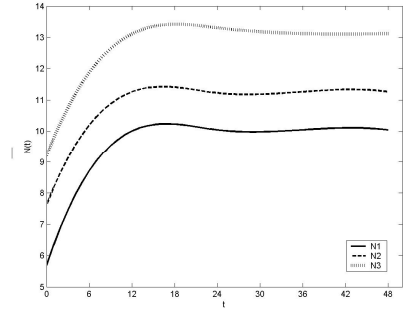}
    \label{fig3b}
  \caption{Change of phase shift. Population size for different values of $\alpha_{K}$
   (N1: $\alpha_{K}=0.1$, N2: $\alpha_{K}=0.5$ and N3: $\alpha_{K}=0.7$).}
    \end{center}
  \end{figure} Often debated questions are the choice of the harvesting strategies and timing of harvesting.
For simplicity, we consider static vs adaptive fishing strategies.\\

{\bf Static Fishing Strategy.}
Consider a fishery manager who has no access to previous fishery data. He starts fishing
(on Fig 6 curve $N1$) all year with an annual quota of 12 tons and realized  that in one-year period  a fishstock is decreased significantly. Then the manager  decides to shorten a fishing season, starts it in June and fishes for six months with a monthly quota of 2 tons (curve $N2$), and thereafter, he starts in September with the quota of 4 tons per month and fishes three months ( curve $N3$ ),
\begin{figure}[ht]
    \begin{center}
    \includegraphics[width=2.5in, viewport=76 303 584 719 clip=true]{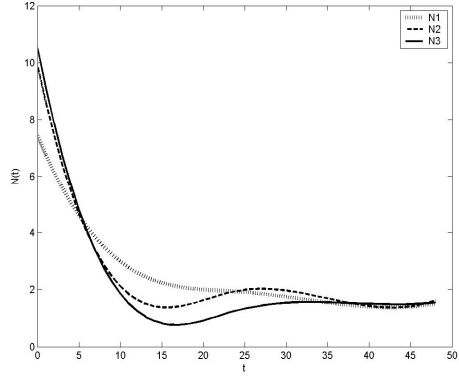}
    \label{fig6}
     \caption{Static Fishing Strategy}
    \end{center}
  \end{figure}
  \\
 \newpage
    {\bf Adaptive Fishing Strategy.}
Let a fishery manager have access (Fig 7 first graph) to the fishery data. He decides to fish in March  because at that time the population attains its maximum. Curve $N2$ represents a six-month fishing, starting in March
with 2 tons per month, whereas on curve $N3$, fishing takes place with the quota of 4 tons per month but in a three-month period. Clearly, the latter strategy is more efficient than the static strategy (curve $N1$ on Fig 7) which
represents a greater risk of depletion of the fishstock.
\begin{figure}[ht]
    \begin{center}
    \includegraphics[width=2.5in, viewport=74 270 530 750,clip=true]{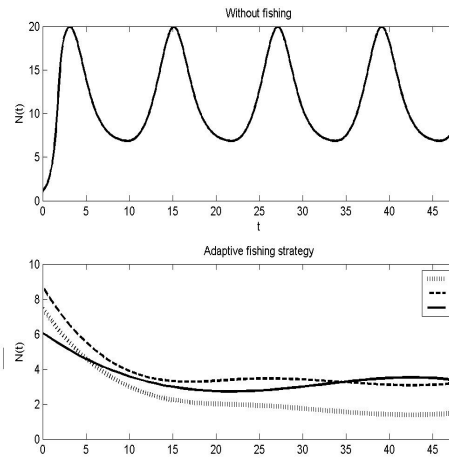}
    \label{fig6}
   \caption{ Adaptive Fishing Strategy}
    \end{center}
  \end{figure}

\section{Discussion}
In this paper we illustrate various effects of the varying
environmental carrying capacity and intrinsic rates on dynamics of
marine populations. Application of the food-limited model with $\gamma$-nonlinearity is consistent with
fishery data, and supports a well-known feature of the population models: average population size decreases
with an increase of amplitude of variation in the carrying capacity $K(t)$, but $\overline{N}\leq \overline{K}$ and the dynamics of the population must ride on this highly variable resource. On other hand,
even large oscillations in $r(t)$ alone leave the system's behavior practically unchanged. Note that similar results were obtained for the Fox fishery model \eqref{Fox} in \cite{Id}. However, the qualitative behavior of the system depends critically on the phase difference between the oscillations in $r(t)$ and $K(t)$. For example, if functions $r(t)$ and $K(t)$ have a $180^{\circ}$ phase shift then
the relationship between the average population size and environmental and demographic variations is
very different. \\
 With no access to fisheries data a static fishing strategy represents a greater risk of extinction of marine populations under severe harvesting, compared with adaptive strategies that save fishstock in a cost effective manner.\\

{\bf Acknowledgements}\\
We wish to express thanks to Dr. D. Barker (Fisheries and Aquaculture department at Vancouver Island University) whose comments significantly improved the text.

\end{document}